\documentclass{amsart}

\usepackage{bbm}
\usepackage{physics} 
\usepackage{mathrsfs} 

\usepackage[utf8]{inputenc} 
\usepackage[T1]{fontenc} 
\usepackage{amsmath} 
\usepackage{amssymb} 
\usepackage{hyperref}


\begin{document}

\title[KMS Inequalities: From Elliptic Operators to Constant Rank]{KMS Inequalities: From Elliptic Operators to Constant Rank}

\author{Paul Stephan}

\address{
    Paul Stephan,
    Universität Konstanz, Universitätsstrasse 10, 78464 Konstanz, Germany
}
\email{paul.stephan@uni-konstanz.de}

\begin{abstract}
Korn-Maxwell-Sobolev (KMS) inequalities represent a tool for estimating differential expressions and have gained particular importance in recent years, especially concerning elliptic operators. In my Master's thesis, together with Peter Lewintan (University of Duisburg-Essen), we extended this concept to also apply to operators of constant rank. This makes it possible to cover more complex structures such as the curl or divergence of vector fields.

A key difference from the elliptic theory is that in the constant rank case, a correction term $\Pi_\mathbb{B}$ is necessary on the left-hand side of the inequality. Results were also obtained for the limiting case $p=1$, although additional assumptions are required here. This article provides an illustrative introduction to KMS inequalities and demonstrates their application in both the elliptic and constant rank cases.
\end{abstract}

\maketitle 

\section{Coercivity in Application}
We first demonstrate the utility and necessity of so-called Korn-Maxwell-Sobolev inequalities using the following application example from the pseudostress formulation of the stationary Stokes equations.
Let $\Omega$ be a domain and $f: \Omega \to \mathbb{R}^3$. We seek a pressure function $\pi: \Omega \to \mathbb{R}$, a velocity $v: \Omega \to \mathbb{R}^3$, and a stress tensor $\sigma: \Omega \to \mathrm{Sym}(3)$ such that the system
\begin{equation}
\begin{cases}
    \sigma - \mu \varepsilon(v) + \pi \operatorname{Id} = 0 \\
    \operatorname{div} \sigma = f \\
    \operatorname{div} v = 0, 
\end{cases} \label{Stokes}
\end{equation}
is satisfied in $\Omega$, where
\[
\varepsilon(v) := \frac{1}{2} \left(\mathrm{D}v + \mathrm{D}v^\top \right).
\]
With the choice
\[
\operatorname{dev} \sigma := \sigma - \frac{1}{3} \operatorname{tr}(\sigma) \operatorname{Id}
\]
the system \eqref{Stokes} can be reformulated as
\begin{equation}
    \begin{cases}
        \operatorname{dev} \sigma - \mu \varepsilon(v) = 0 \\
        \pi = - \frac{\operatorname{tr} \sigma}{3} \\
        \operatorname{div} \sigma = f,
    \end{cases}
\end{equation}
where the divergence-free condition of $v$ is now reflected in the trace-free nature of the first equation.
Since $\pi$ is thus completely determined by $\sigma$, the second line can be disregarded for a moment, and we can otherwise describe the problem variationally. This means we seek minimizers of the functional
\[
\mathcal{F}[v; \Omega] = \int \limits_\Omega \| \operatorname{dev} \sigma - \mu \varepsilon(v) \|^2 \dd{x} + \int \limits_\Omega \|\operatorname{div} \sigma - f\|^2 \dd{x}
\]

on suitable function spaces.
To apply established techniques to this minimization problem, such as the direct method of the calculus of variations, \emph{coercivity} is a crucial property. This means we are interested in applying an inequality of the form
\[
| \mathcal{F}[v; \Omega] | \geq C \| v \|
\]
for some $C>0$. Such an inequality can be formulated in the present case as follows:
 \begin{align}
\int \limits_\Omega \| \operatorname{dev} \sigma - \mu \varepsilon(v) \|^2 \dd{x} & + \int \limits_\Omega \|\operatorname{div} \sigma - f\|^2 \dd{x} \label{eq:koerz}\\
&\geq C \left(\| \operatorname{div} v \|_{Z_1}^2 + \|v\|_{Z_2}^2 \right), \notag
 \end{align}
 where $Z_1, Z_2$ are suitable $\mathrm{L}^p$-valued spaces.
Proving such inequalities is generally possible, as demonstrated in \cite[Section 7.2]{Models}. However, this involves a certain effort. The question in the context of Korn–Maxwell–Sobolev inequalities, which were first mentioned under this name in \cite{gmeineder2021korn}, concerns the possibility of generalizing such inequalities to apply them to the broadest possible class of differential operators. Furthermore, Korn-Maxwell-Sobolev inequalities make it possible to obtain the inequality \eqref{eq:koerz} by verifying simple algebraic conditions instead of performing concrete calculations. They also find applications, for example, in micromorphic models by \textsc{Neff} et al.
(see \cite{Neff_Relaxed2, Neff_Relaxed_2015, Neff_Unifying_2014, Sky_Primal_2022}),
in gradient plasticity by \textsc{Neff} et al. (see \cite{plastic_gradient1, plastic_gradient2, plastic_gradient3}),
in incompatible elasticity by \textsc{Amstutz \& Van Goethem}
(see \cite{Amstutz_Incompatibility_2017, Amstutz_Existence_2020}), or in planar elasticity by \textsc{Arnold} et al.
(see \cite{Arnold_Family_1984}), to name just a few. For an overview of further applications, see also \cite{Lewintan_Korn_2021}.

\section{Homogeneous Differential Operators}
Let $\text{Lin}(\mathbb{R}^d; \mathbb{R}^l)$ be the space of all linear maps from $\mathbb{R}^d$ to $\mathbb{R}^l$. Given a family $\{B_\alpha\}_{|\alpha| = k}$ of linear maps $B_{\alpha} \in \text{Lin}(\mathbb{R}^d; \mathbb{R}^l)$, indexed by multi-indices $\alpha$ with $|\alpha| = k$. For a function $v: \mathbb{R}^n \to \mathbb{R}^d$, we define the homogeneous operator $\mathbb{B}$ as
\[
\mathbb{B}v := \sum_{|\alpha| = k} B_{\alpha} \, \partial^\alpha v.
\]
The corresponding Fourier symbol map is given by
\[
\mathbb{B}[\xi] := \sum_{|\alpha| = k} B_{\alpha} \, \xi^\alpha,
\]
where $\xi^\alpha = \xi_1^{\alpha_1} \cdots \xi_n^{\alpha_n}$ for $\xi = (\xi_1, \ldots, \xi_n)$ and $\alpha = (\alpha_1, \ldots, \alpha_n)$.

The operator $\mathbb{B}$ is called elliptic if for all $\xi \in \mathbb{S}^{n-1}$ the map
\[
\mathbb{B}[\xi] : \mathbb{R}^d \to \mathbb{R}^l
\]
is injective. $\mathbb{B}$ is said to have constant rank if there exists an $r \in \mathbb{N}$ such that
\[
\operatorname{rank} \mathbb{B}[\xi] = r
\]
holds for all $\xi \in \mathbb{S}^{n-1}$. \\
We now wish to focus on investigating the beneficial properties of elliptic operators in the context of coercive estimates. In particular, we want to present Korn's inequality for elliptic operators, which plays a central role. The generalized Korn inequality can be stated as follows:
\begin{equation}
    \|\mathrm{D} P\|_{\mathrm{L}^p(\mathbb{R}^n; \mathbb{R}^{d \times n})} \lesssim \|\mathbb{B} P\|_{\mathrm{L}^p(\mathbb{R}^n; \mathbb{R}^l)}, \label{Korn:ell}
\end{equation}
which holds for all test functions \( P \in \mathrm{C}_c^\infty(\mathbb{R}^n; \mathbb{R}^d) \) and exponents $1 < p < \infty$. Here and in the following, ``$\lesssim$'' expresses that the inequality holds up to a constant factor $C>0$, which depends only on $n, p, d, l, k$ and $\mathbb{B}$ (but not on $P$). \\
This inequality demonstrates that the derivative $\mathrm{D}P$ of a function $P$ can be controlled in the $\mathrm{L}^p$-norm by the action of an elliptic differential operator $\mathbb{B}$ on $P$, at least if $p>1$. We will return to the case $p=1$ in the last section. \\
The derivation of \eqref{Korn:ell} is based on Mihlin's theorem (see \cite{Duoandikoetxea2001, Stein1970}): Let $m \in \mathrm{C}^{n+2}(\mathbb{R}^n \setminus \{0\})$ be a symbol that satisfies the condition
\[
    |\partial^\alpha m(\xi)| \lesssim |\xi|^{-|\alpha|}
\]
for all multi-indices $\alpha$ with $|\alpha| \leq n+2$. Then Mihlin's theorem states that the associated Fourier multiplier $T_m$ is bounded as an operator on $\mathrm{L}^p(\mathbb{R}^n)$ for $1<p<\infty$. This theorem provides the necessary tool to derive Korn's inequality via Fourier methods. \\
The basic idea of this derivation can be illustrated in the whole space by representing the derivative of $P$ with
\[
    \partial^\alpha P = \mathcal{F}^{-1} \left( \underset{{=:m(\xi)}}{\underbrace{ (i\xi)^\alpha \left( \mathbb{B}^*[\xi] \mathbb{B}[\xi] \right)^{-1} {\mathbb{B}[\xi]}^*}} \widehat{\mathbb{B} P}(\xi) \right) 
\]
where the choice of the multiplier $m(\xi)$ justifies the use of Mihlin's theorem and yields the statement of Korn's inequality \eqref{Korn:ell}. 

\section{From Korn to Korn-Maxwell-Sobolev}
If we choose the operator $\mathbb{B}=\varepsilon$ in \eqref{Korn:ell} with
\[
\varepsilon(u):= \frac{1}{2} \left( \mathrm{D}u + \mathrm{D}u^\top \right)
\]
as the symmetric gradient, we obtain the classical form of Korn's inequality, which is of central importance in elasticity theory.
For $u \in \mathrm{C}_c^\infty(\mathbb{R}^n; \mathbb{R}^n)$ and $1 < p < \infty$, it then holds that:
\begin{equation}
\| \mathrm{D} u \|_{\mathrm{L}^p(\mathbb{R}^n; \mathbb{R}^{n \times n})} \lesssim \| \varepsilon(u) \|_{\mathrm{L}^p(\mathbb{R}^n; \mathbb{R}^{n \times n})}. \label{Korn}
\end{equation}
This inequality is surprising because here $n^2$ entries on the left-hand side can be controlled by $\frac{n(n+1)}{2}$ entries on the right-hand side. It plays a fundamental role because, in particular, numerous mechanical phenomena depend solely on the symmetric part of the gradient, especially for linearly elastic materials. The reason for this is that the response of an elastic body is typically modeled only by the symmetric part of the gradient of the displacement field. The antisymmetric part, on the other hand, causes no stretching or compression in the material but merely describes an infinitesimal rigid rotation. This can be expressed by decomposing the derivative $\mathrm{D} u$ of a displacement field $u$ into symmetric and antisymmetric parts:
\[
\mathrm{D} u = \varepsilon(u) + \omega(u)
\]
where
\[
\omega(u) = \frac{1}{2}(\mathrm{D} u - \mathrm{D} u^\top)
\]
denotes the antisymmetric part. The symmetric part $\varepsilon(u)$ describes the actual deformation of the material, which leads to the storage of elastic energy.
 \\
Subsequently, the question arises whether a similar inequality to \eqref{Korn} is possible for general vector fields $P$, where $\varepsilon$ is the projection onto the symmetric part:
\[
\| P \|_{\mathrm{L}^p(\mathbb{R}^n; \mathbb{R}^{n \times n})} \lesssim \| \operatorname{sym}[P] \|_{\mathrm{L}^p(\mathbb{R}^n; \mathbb{R}^{n \times n})},
\]
for all $P \in \mathrm{C}_c^\infty(\mathbb{R}^n; \mathbb{R}^{n \times n})$. This is not the case. To illustrate, let $P$ be a skew-symmetric matrix. In this case, the symmetric part $\operatorname{sym}[P] = 0$, so the right-hand side of the inequality vanishes, while the left-hand side is not necessarily zero. This shows that the norms of general fields cannot be controlled by the symmetric part alone, which means there must be a special property of $P = \mathrm{D}u$ (like being a gradient) that plays a role in this context. \\
To circumvent the above limitation, for $1<p<n$, there is the possibility of choosing a new approach:
\begin{align}
& \| P \|_{\mathrm{L}^{p^*}(\mathbb{R}^n; \mathbb{R}^{n \times n})} \label{KMS:sym} \\
& \qquad \lesssim \| \operatorname{sym}[P] \|_{\mathrm{L}^{p^*}(\mathbb{R}^n; \mathbb{R}^{n \times n})} + \| \operatorname{Curl}(P) \|_{\mathrm{L}^p(\mathbb{R}^n; \mathbb{R}^l)}, \notag
\end{align}
where for scaling reasons $p^*=\frac{np}{n-p}$ must hold (the Sobolev conjugate exponent).
The approach described here leads to the development of Korn-Maxwell-Sobolev inequalities, which allow consideration of both the symmetric part and the curl of tensor fields, thus covering a broad class of operators. Based on results by \textsc{Neff, Pauly}, and \textsc{Witsch} (see \cite{neff_pauly1, neff_pauly2}),
\textsc{Gmeineder, Lewintan}, and \textsc{Neff}
 were able to prove in \cite{AnalyticSplit} that for $n \ge 3$, $1 \leq p < n$, and a linear map $\mathscr{A} : \mathbb{R}^{m \times n} \to \mathbb{R}^N$ such that the corresponding homogeneous differential operator $\mathbb{A}$ (acting on $m \times n$ tensor fields) is elliptic, the following estimate holds:
\begin{align}
& \| P \|_{\mathrm{L}^{p^*}(\mathbb{R}^n; \mathbb{R}^{m \times n})} \label{ASplit} \\
& \qquad \lesssim \| \mathscr{A}[P] \|_{\mathrm{L}^{p^*}(\mathbb{R}^n; \mathbb{R}^N)} + \| \operatorname{Curl}(P) \|_{\mathrm{L}^p(\mathbb{R}^n; \mathbb{R}^l)}, \notag
\end{align}
for all $P \in \mathrm{C}^\infty_c(\mathbb{R}^n; \mathbb{R}^{m \times n})$. The converse is also valid, meaning that if \eqref{ASplit} holds for all $P \in \mathrm{C}_c^\infty(\mathbb{R}^n; \mathbb{R}^{m \times n})$, then $\mathbb{A}$ must already be elliptic. Here, they generalized a result by \textsc{Gmeineder} and \textsc{Spector} from \cite{gmeineder2021korn}. \\
Inequality \eqref{KMS:sym} appears as a special case of inequality \eqref{ASplit} for the case $\mathscr{A}=\operatorname{sym}$ (with $m=n$). Both inequalities fall into the category of Korn-Maxwell-Sobolev inequalities, whose name we want to explain in more detail: If we set $P = \mathrm{D}u$ for a vector field $u \in \mathrm{C}_c^\infty(\mathbb{R}^n; \mathbb{R}^n)$ and define $\mathscr{A}[P] = \frac{1}{2}(P + P^\top) = \varepsilon(u)$, we recover the classical Korn inequality \eqref{Korn} that we saw above (or \eqref{Korn:ell} for an arbitrary elliptic $\mathbb{A}$ acting on gradients, noting $\mathrm{Curl}(\mathrm{D}u)=0$).
The Maxwell component refers to the frequently used \emph{Curl} operator in Maxwell's equations. An example of this is the induction equation
\[
\operatorname{Curl} E = - \frac{\partial B}{\partial t},
\]
which describes how a time-varying magnetic field $B$ induces an electric field $E$.
Finally, "Sobolev" refers to the exponent $p^*$, which appears on the left-hand side of the inequality as the Sobolev exponent. If one sets $\mathscr{A} = 0$, the inequality reduces to a Sobolev-type inequality.

\section{An Extended Version of KMS Inequalities}

Having considered an inequality with the \emph{Curl} operator on the right-hand side, the question now arises whether KMS inequalities can also be generalized for a larger class of operators. Indeed, \textsc{Gmeineder, Lewintan}, and \textsc{Neff} achieved a crucial result here in \cite{article}, where they could replace the Curl operator with any homogeneous operator $\mathbb{B}$ that is elliptic on a suitable set: \\
Let $\mathscr{A} : \mathbb{R}^d \to \mathbb{R}^N$ be a linear map and $\mathbb{B}$ be a homogeneous differential operator of order $k$, which is elliptic on $\ker(\mathscr{A})$. Then for $1 < p < n$ and $P \in \mathrm{C}_c^\infty(\mathbb{R}^n; \mathbb{R}^d)$, the estimate holds:
\begin{align}
& \| P \|_{\dot{\mathrm{W}}^{k-1, p^*}(\mathbb{R}^n; \mathbb{R}^d)} \label{Korn:ellip} \\ 
& \qquad \lesssim \| \mathscr{A}[P] \|_{\dot{\mathrm{W}}^{k-1, p^*}(\mathbb{R}^n; \mathbb{R}^N)} + \| \mathbb{B}P \|_{\mathrm{L}^p(\mathbb{R}^n; \mathbb{R}^l)}. \notag
\end{align}
Remarkably, the converse also holds: If \eqref{Korn:ellip} is satisfied, then $\mathbb{B}$ must be elliptic on the kernel of $\mathscr{A}$. \\
Since the Curl operator itself is not elliptic (it has a non-trivial kernel), it is natural to try to consider an even more general class of differential operators on the right-hand side, namely those of constant rank. A key challenge to overcome here is the following: If the previous inequality is already equivalent to the ellipticity of $\mathbb{B}$ on the kernel of $\mathscr{A}$, then the inequality in the constant rank case cannot look \emph{exactly} the same. Modifications are therefore necessary.

\section{The Main Result: KMS Inequalities for Operators with Constant Rank}

Now the necessary theoretical foundations have been established, so that we can present our main result in its full generality. \\
Let $\mathscr{A} : \mathbb{R}^d \to \mathbb{R}^N$ be a linear map and $\mathbb{B}$ be a homogeneous differential operator of order $k$, which has constant rank on $\ker(\mathscr{A})$. Then for $1 < p < n$ and $P \in \mathrm{C}_c^\infty(\mathbb{R}^n; \mathbb{R}^d)$, the estimate holds:
\begin{align}
& \| P - \Pi_{\mathbb{B}}\mathcal{P}_{\ker(\mathscr{A})}[P] \|_{\dot{\mathrm{W}}^{k-1, p^*}(\mathbb{R}^n; \mathbb{R}^d)} \label{Korn:const} \\ 
& \qquad \lesssim
\| \mathscr{A}[P] \|_{\dot{\mathrm{W}}^{k-1, p^*}(\mathbb{R}^n; \mathbb{R}^N)} + \| \mathbb{B}P \|_{\mathrm{L}^p(\mathbb{R}^n; \mathbb{R}^l)} \notag,
\end{align}

where $\Pi_{\mathbb{B}}: \mathrm{L}^p \to \mathrm{L}^p$ denotes a projection associated with $\mathbb{B}$ and $\mathcal{P}_{\ker(\mathscr{A})}$ denotes the projection onto the kernel of $\mathscr{A}$. As can be seen, the right-hand side of \eqref{Korn:const} does not change compared to \eqref{Korn:ellip}, but instead, the correction term is newly added on the left-hand side. \\
One might now ask whether, if \eqref{Korn:const} is valid for all $P \in \mathrm{C}_c^\infty(\mathbb{R}^n; \mathbb{R}^d)$, the operator $\mathbb{B}$ must already have constant rank on the kernel of $\mathscr{A}$. We could not yet definitively clarify this question, but instead, we were able to show that the inequality
\begin{align}
& \| P - \Pi_{\mathbb{B}}\mathcal{P}_{\ker(\mathscr{A})}[P] \|_{\dot{\mathrm{W}}^{-1, p^*}(\mathbb{R}^n; \mathbb{R}^d)} \label{Korn:const2} \\
& \qquad \lesssim
\| \mathscr{A}[P] \|_{\dot{\mathrm{W}}^{-1, p^*}(\mathbb{R}^n; \mathbb{R}^N)} + \| \mathbb{B}P \|_{\dot{\mathrm{W}}^{-k,p}(\mathbb{R}^n; \mathbb{R}^l)} \notag,
\end{align}
for all $P \in \mathrm{C}_c^\infty(\mathbb{R}^n; \mathbb{R}^d)$ is equivalent to $\mathbb{B}$ having constant rank on $\ker(\mathscr{A})$. Here, the backward direction follows directly from inequality \eqref{Korn:const}, while the forward direction builds on a result by \textsc{Guerra} and \textsc{Rai\c{t}\u{a}} from \cite{GR}. \\
A special case occurs when $\ker(\mathscr{A}) \subset \ker(\mathbb{B}[\xi])$ for all $\xi \neq 0$. This is the case, for example, in three dimensions when $\mathscr{A} = \operatorname{dev}$ and $\mathbb{B}=\mathrm{sym} \mathrm{Curl}$, since here $\mathcal{P}_{\ker(\mathscr{A})} P = \frac{\operatorname{tr} P}{3} \mathbbm{1}_3$ and thus
\[
\operatorname{dev} P = P - \frac{\operatorname{tr} P}{3} \mathbbm{1}_3 = \mathcal{P}_{(\ker(\mathscr{A}))^\perp } P
\]
holds (see \cite[Remark 3.2]{preprint}). In this case, \eqref{Korn:const} reduces to:
\begin{align}
& \| \mathcal{P}_{(\ker(\mathscr{A}))^\perp }[P] \|_{\dot{\mathrm{W}}^{k-1, p^*}(\mathbb{R}^n; \mathbb{R}^d)} \notag \\ 
& \qquad \lesssim
\| \mathscr{A}[P] \|_{\dot{\mathrm{W}}^{k-1, p^*}(\mathbb{R}^n; \mathbb{R}^N)} + \| \mathbb{B}P \|_{\mathrm{L}^p(\mathbb{R}^n; \mathbb{R}^l)} \notag.
\end{align}

\section{Sketch of the Proof of the Main Result}

The basis of our proof is a result by \textsc{Fonseca} and \textsc{Müller} from \cite{FM}. This states that for a homogeneous differential operator $\mathbb{B}$ of order $k$ with constant rank, there exists a constant $C_p > 0$ such that
\begin{equation}
\| P - \Pi_\mathbb{B} P \|_{\mathrm{L}^p(\mathbb{R}^n; \mathbb{R}^d)} \leq C_p \|\mathbb{B} P\|_{\mathrm{W}^{-k,p}(\mathbb{R}^n; \mathbb{R}^l)} \label{res:FM}
\end{equation}
is satisfied for all $P \in \mathrm{C}_c^\infty(\mathbb{R}^n; \mathbb{R}^d)$.
The idea now is to split the test function $P$ into a part on the kernel of $\mathscr{A}$ and a part on the complement of this set. More precisely, we write
\[
P = \mathcal{P}_{\ker(\mathscr{A})}[P] + \mathcal{P}_{\ker(\mathscr{A})^\perp}[P]
\]
and estimate the two terms separately. For the second term ($\mathcal{P}_{\ker(\mathscr{A})^\perp}[P]$), we use that $\mathscr{A}$ is injective on $\ker(\mathscr{A})^\perp$, so that pointwise
\begin{equation}
|\mathcal{P}_{\ker(\mathscr{A})^\perp}[P]| \lesssim |\mathscr{A}[\mathcal{P}_{\ker(\mathscr{A})^\perp}[P]]| = |\mathscr{A}[P]| \label{absch:pkt}
\end{equation}
holds (since $\mathscr{A}[\mathcal{P}_{\ker(\mathscr{A})}[P]]=0$). The injectivity is seen as follows: If there are two elements $\xi_1, \xi_2 \in \ker(\mathscr{A})^\perp$ with $\mathscr{A}[\xi_1] = \mathscr{A}[\xi_2],$ then $\mathscr{A}[\xi_1 - \xi_2]= 0$ and thus $\xi_1 - \xi_2 \in \ker(\mathscr{A}) \cap \ker(\mathscr{A})^\perp = \{0\}.$
For the other term ($\mathcal{P}_{\ker(\mathscr{A})}[P]$), we can apply the constant rank theory to the operator $\mathbb{B}$ restricted to functions taking values in $\ker(\mathscr{A})$. Let $P_{\ker} := \mathcal{P}_{\ker(\mathscr{A})}[P]$. Since $\mathbb{B}$ has constant rank on $\ker(\mathscr{A})$, we can apply a version of \eqref{res:FM} to $P_{\ker}$:
\begin{align}
\| P_{\ker} - & \Pi_{\mathbb{B}} P_{\ker} \|_{\mathrm{L}^p(\mathbb{R}^n; \mathbb{R}^d)} \\
& \lesssim \| \mathbb{B} P_{\ker} \|_{\dot{\mathrm{W}}^{-k,p}(\mathbb{R}^n; \mathbb{R}^l)} \notag \\
    & = \| \mathbb{B}(P-\mathcal{P}_{\ker(\mathscr{A})^\perp}[P]) \|_{\dot{\mathrm{W}}^{-k,p}(\mathbb{R}^n; \mathbb{R}^l)} \notag \\
    & \lesssim \| \mathbb{B}P \|_{\dot{\mathrm{W}}^{-k,p}(\mathbb{R}^n; \mathbb{R}^l)} + \| \mathbb{B} \mathcal{P}_{\ker(\mathscr{A})^\perp}[P] \|_{\dot{\mathrm{W}}^{-k,p}(\mathbb{R}^n; \mathbb{R}^l)} \notag \\
    & \lesssim \| \mathbb{B}P \|_{\dot{\mathrm{W}}^{-k,p}(\mathbb{R}^n; \mathbb{R}^l)} + \| \mathcal{P}_{\ker(\mathscr{A})^\perp} [P]\|_{\mathrm{L}^p(\mathbb{R}^n; \mathbb{R}^d)} \notag \quad (\text{since } \mathbb{B} \text{ is order } k) \\
    & \lesssim \| \mathbb{B}P \|_{\dot{\mathrm{W}}^{-k,p}(\mathbb{R}^n; \mathbb{R}^d)} + \| \mathscr{A}[P]\|_{\mathrm{L}^p(\mathbb{R}^n; \mathbb{R}^l)}, \notag
\end{align}
where in the last step we used \eqref{absch:pkt} integrated in $\mathrm{L}^p$. \\
If we insert the $(k-1)$-th derivatives of $P$ into the above inequalities, we obtain 
\begin{align*}
& \| (P - \Pi_\mathbb{B} \mathcal{P}_{\ker(\mathscr{A})}[P]) \|_{\dot{\mathrm{W}}^{k-1, p^*}(\mathbb{R}^n; \mathbb{R}^d)} \\ 
& \qquad \lesssim \| \mathscr{A}[P] \|_{\dot{\mathrm{W}}^{k-1, p^*}(\mathbb{R}^n; \mathbb{R}^N)} + \| \mathbb{B}  P \|_{\dot{\mathrm{W}}^{-1, p^*}(\mathbb{R}^n; \mathbb{R}^l)}.
\end{align*} Now it suffices to note that by classical Sobolev embedding for $1<q<n$, $\dot{\mathrm{W}}^{1,q}(\mathbb{R}^n; \mathbb{R}^l)$ embeds into $\mathrm{L}^{\frac{nq}{n-q}}(\mathbb{R}^n; \mathbb{R}^l)$, which by duality implies that for
\[
q = \frac{np}{np-n+p}
\]
the space $\mathrm{L}^p(\mathbb{R}^n; \mathbb{R}^l)$ embeds into the space
\[ 
\dot{\mathrm{W}}^{-1, \frac{q}{q-1}}(\mathbb{R}^n; \mathbb{R}^l),
\]
where
\[
\frac{q}{q-1} = \frac{np}{n-p}
\]
holds.
It follows:
\begin{align*}
& \| P - \Pi_\mathbb{B} \mathcal{P}_{\ker(\mathscr{A})}[P] \|_{\dot{\mathrm{W}}^{k-1, p^*}(\mathbb{R}^n; \mathbb{R}^d)} \\ 
& \qquad \lesssim \| \mathscr{A}[P] \|_{\dot{\mathrm{W}}^{k-1, p^*}(\mathbb{R}^n; \mathbb{R}^N)} + \| \mathbb{B} P \|_{\mathrm{L}^p(\mathbb{R}^n; \mathbb{R}^l)}
\end{align*}
and thus the claimed inequality.

\section{The Case \texorpdfstring{$p=1$}{p=1}}
For the case $p=1$, the above argument cannot be directly transferred, as the classical theory of Fourier multipliers and Mihlin's theorem is not available here. Indeed, it is known from Ornstein's non-inequality (see \cite{ornstein1962non}) that counterexamples can even be constructed in this case and that further conditions must be imposed on $\mathbb{B}$. \\
An important condition is that of being a \emph{cancelling operator}, which goes back to \textsc{Van Schaftingen} (see \cite{vanschaftingen2013limiting}). Here, we call a differential operator $\mathbb{B}$ \emph{cancelling} if and only if
\[
\bigcap \limits_{\xi \neq 0} \operatorname{Im}(\mathbb{B}[\xi]) = \{0\} 
\]
holds. This definition is central because the inequality
\[
\|P\|_{\mathrm{W}^{k-1, \frac{n}{n-1}}(\mathbb{R}^n;\mathbb{R}^d)} \lesssim \|\mathbb{B} P\|_{\mathrm{L}^1(\mathbb{R}^n; \mathbb{R}^l)}
\]
for all $P \in \mathrm{C}_c^\infty(\mathbb{R}^n; \mathbb{R}^d)$ is equivalent to $\mathbb{B}$ being elliptic and cancelling (see \cite[Theorem 1.3]{vanschaftingen2013limiting}). This means, in particular, that for this form of inequality in the case $p=1$, ellipticity is no longer sufficient, which implies that this is also the case for KMS inequalities. Indeed, \textsc{Gmeineder, Lewintan}, and \textsc{Van Schaftingen} were able to show in \cite{GLVS} that the inequality
\begin{align}
& \| P \|_{\dot{\mathrm{W}}^{k-1, \frac{n}{n-1}}(\mathbb{R}^n; \mathbb{R}^d)} \label{Korn:ellipP=1} \\
& \qquad \lesssim \| \mathscr{A}[P] \|_{\dot{\mathrm{W}}^{k-1, \frac{n}{n-1}}(\mathbb{R}^n; \mathbb{R}^N)} + \| \mathbb{B}P \|_{\mathrm{L}^1(\mathbb{R}^n; \mathbb{R}^l)} \notag
\end{align}
is equivalent to $\mathbb{B}$ being elliptic and cancelling on the kernel of $\mathscr{A}$. This is therefore an optimal extension of \eqref{Korn:ellip} to the case $p=1$. \\
Here too, it was the task of my Master's thesis to transfer this to the case $p=1$. The investigation showed that under the assumption that $\mathbb{B}$ has constant rank on the kernel of $\mathscr{A}$ and is cancelling there, it holds that
\begin{align}
& \| P - \Pi_{\mathbb{B}}\mathcal{P}_{\ker(\mathscr{A})}[P] \|_{\dot{\mathrm{W}}^{k-1, \frac{n}{n-1}}(\mathbb{R}^n; \mathbb{R}^d)} \label{Korn:constP=1} \\
& \qquad \lesssim
\| \mathscr{A}[P] \|_{\dot{\mathrm{W}}^{k-1, \frac{n}{n-1}}(\mathbb{R}^n; \mathbb{R}^N)} + \| \mathbb{B}P \|_{\mathrm{L}^1(\mathbb{R}^n; \mathbb{R}^l)} \notag,
\end{align}
for all $P \in \mathrm{C}_c^\infty(\mathbb{R}^n; \mathbb{R}^d)$, which thus also transfers \eqref{Korn:const} to the case $p=1$. The question of whether the conditions that $\mathbb{B}$ must have constant rank on the kernel of $\mathscr{A}$ and be cancelling there are also \emph{necessary} is currently still open. \\
Although these prerequisites may seem very abstract at first glance, they can certainly be illustrated using concrete examples. If $\mathbb{B} = \operatorname{Curl}$, it is known that $\mathbb{B}$ has global constant rank and is globally cancelling. Consequently, \eqref{Korn:constP=1} holds for any linear map $\mathscr{A}$, for example $\mathscr{A} = \operatorname{tr}$. In this case, $\Pi_\mathbb{B}$ can be calculated explicitly and we get:
\[
\Pi_\mathbb{B} \mathcal{P}_{\ker(\mathscr{A})} P (x) = \frac{1}{n \omega_n} \int \limits_{\mathbb{R}^n} \frac{x-y}{|x-y|^n} \otimes (\operatorname{Div} \operatorname{dev} P)(y) \dd{y} 
\]
or in view of \eqref{Korn:constP=1} (with $k=1$)
\begin{align*}
     &\left( \int \limits_{\mathbb{R}^n} \left| P(x)-   \frac{1}{n \omega_n} \int \limits_{\mathbb{R}^n} \frac{x-y}{|x-y|^n} \otimes (\operatorname{Div} \operatorname{dev} P)(y) \dd{y} \right|^{\frac{n}{n-1}} \dd{x} \right)^\frac{n-1}{n} \\
     & \quad \lesssim \left( \left( \int \limits_{\mathbb{R}^n} |\operatorname{tr} P(x)|^\frac{n}{n-1} \dd{x} \right)^\frac{n-1}{n} + \int \limits_{\mathbb{R}^n} |\operatorname{Curl} P(x) | \dd{x} \right).
\end{align*}

\section{Outlook - KMS Inequalities on Domains}
After having formulated and treated Korn-Maxwell-Sobolev inequalities only on the whole space so far, we want to address the question in this last section of what can be said on domains. Meanwhile, there are initial results concerning Korn inequalities on domains, which are equivalent to $\mathbb{A}$ being $\mathbb{C}$-elliptic, see \cite[Theorem 1.3]{Diening_Sharp_2024}. Here, we call $\mathbb{A}$ $\mathbb{C}$-elliptic if
$\mathbb{A}[\xi]$ is injective for all $\xi \in \mathbb{C}^n \setminus \{0\}$. This condition is thus more restrictive than the previously assumed ellipticity, where we only required $\xi \in \mathbb{R}^n \setminus \{0\}$. Sobolev inequalities also exist on domains (see \cite{Gmeineder_Embeddings_2019}). It would be an interesting question to what extent these results can also be transferred to KMS inequalities.

 \bibliographystyle{amsplain}
 \bibliography{ArtikelDMV}

\end{document}